\documentclass[a4paper,10pt]{article}
\usepackage{mathtext}
\usepackage[T1,T2A]{fontenc}
\usepackage[english]{babel}
\usepackage[OT4]{fontenc}
\usepackage[cp1250]{inputenc}
\usepackage{amsmath}
\usepackage{amsfonts}
\usepackage{amssymb}
\usepackage{mathrsfs}
\usepackage{amsthm}
\usepackage{euscript}
\usepackage{bbm}
\usepackage{graphicx}
\usepackage{bm}
\usepackage{authblk}

\DeclareMathOperator{\supp}{supp}

\newcommand{\degz}[1]{deg_z\left({#1}\right)}

\renewcommand{\leq}{\leqslant}
\renewcommand{\geq}{\geqslant}

\theoremstyle{plain}\newtheorem{T1}{Theorem}
\theoremstyle{plain}
\theoremstyle{plain}\newtheorem{L1}{Lemma}
\theoremstyle{plain}
\theoremstyle{plain}\newtheorem{Def}{Definition}
\theoremstyle{plain}
\theoremstyle{plain}

\author[1]{Krystian Kazaniecki\thanks{Krystian.Kazaniecki@mimuw.edu.pl}}
\author[2]{Michal Wojciechowski\thanks{M.Wojciechowski@impan.pl}}
\affil[1]{Institute of Mathematics, University of Warsaw}
\affil[2]{Institute of Mathematics, Polish Academy of Sciences}

\title{ On the equivalence between the sets of the
trigonometric polynomials}

\begin{document}
\maketitle
\begin{abstract}
In this paper we construct an injection from the linear space of trigonometric polynomials defined on $\mathbb{T}^d$ with bounded degrees with respect to each variable to a suitable linear subspace $L^1_E\subset L^1(\mathbb{T})$. We give such a quantitative condition on $L^1_E$ that this injection is a isomorphism of a Banach spaces equipped with $L^1$ norm and the norm of the isomorphism is independent on the dimension $d$.  
\end{abstract}
\subsection*{Introduction}
The purpose of this article is to study the equivalence between sets of trigonometric polynomials defined on $\mathbb{T}$ with sets of trigonometric polynomials in higher dimensions. For a given vector of integers \mbox{$\tau=(\tau_1, \ldots, \tau_d)$} we want to study operator
\begin{equation*}
f(x)=\sum_{\lambda} a_{\lambda} e^{2\pi i \langle \lambda, x\rangle} \rightarrow Tf(z)= \sum_{\lambda} a_{\lambda} e^{2\pi i \langle \lambda, \tau\rangle z} \qquad \forall x\in\mathbb{T}^d\;\forall z\in\mathbb{T}
\end{equation*}
with some bounds on vectors $\lambda=(\lambda_1,\ldots,\lambda_d)$. We want to find a quantitative criterion on $\tau$ for $T$ being isomorphism in $L^1$ norm. To give a precise formulation we introduce following notation.
\begin{Def} As a $L^p_A(\mathbb{T}^k)$ we will denote a subspace of Banach space $L^p(\mathbb{T}^k)$ defined below
\[L^p_A(\mathbb{T}^k)=\{f\in L^p(\mathbb{T}^k): \supp \widehat{f} \subset A\} \]
\end{Def}
\begin{Def} For a given sequence of natural numbers $\left(a_n\right)_{n\in\mathbb{N}}$ and a sequence of integers $\left(\tau_n\right)_{n\in\mathbb{N}}$ we define sets $E\subset \mathbb{Z}$ and $F\subset\mathbb{Z}^{\mathbb{N}}$, where by $ \mathbb{Z}^{\mathbb{N}}$ we denote a dual group 
to $\mathbb{T}^{\mathbb{N}}$, in the following way:
\begin{equation}
 \begin{split}
 F&=\{\bm{\lambda}\in\mathbb{Z}^{\mathbb{N}}\; :\; |\lambda_n|\leq |a_n|\}\\,
 E&=\{\beta\in\mathbb{Z}\; :\; \; \beta =\sum_{k=1} \tau_k\lambda_k \mbox{ for } \lambda_n\in F\}.
\end{split}
\end{equation}
\end{Def}
Using that notation we can state the main Theorem of this article :
\begin{T1}
 For a given sequence of natural numbers $(a_n)_{n\in\mathbb{N}}$ and a sequence of integers $\tau_n$ satisfying 
 \begin{equation}\label{wlasnosci}
 \begin{split}
  \tau_{k+1}\geq 3 a_{k} |\tau_{k}|&\qquad\forall k\in\mathbb{N},\\
  \sum_{j=1}^{\infty} \frac{|a_j||\tau_j||a_{j+1}|}{|\tau_{j+1}|}&<\infty.
\end{split}
  \end{equation}
Then operator $T : L^1_F(\mathbb{T}^{\mathbb{N}}) \rightarrow L^1_{E}(\mathbb{T})$ given by the formula
\begin{equation}
Tf(x)= \sum_{\bm{\lambda}\in F} \widehat{f}(\bm{\lambda}) e^{2\pi i \langle \bm{\lambda}, \tau \rangle x}
\end{equation}
is an isomorphism, moreover
\[ K^{-1}\|f\|_{L^1_F(\mathbb{T}^{\mathbb{N}})}\leq \| Tf\|_{L^1_E(\mathbb{T})}\leq K\|f\|_{L^1_F(\mathbb{T}^{\mathbb{N}})}\]
with the constant $K$ depending only on the value of \;$\sum_{j=1}^{\infty} \frac{|a_j||\tau_j||a_{j+1}|}{|\tau_{j+1}|}$
\end{T1}
Similar criterion was used by Y. Meyer in case of $L^{p}$ norm and a sequence $a_k\equiv 1$ (cf. \cite{MR0240563}) and later stronger condition was obtained by M. D{\'e}champs (cf. \cite{MR641858}) . However in the case of $L^p$ norm with $1\leq p<2$ both proofs seem to be incomplete. We fix this problem with an elementary proof for $L^1$ norm. For examples of use of such a criterion one can check \cite{info},\cite{MR1649869}.
\\In order to prove the main Theorem we will prove three lemmas. 
\subsection*{Auxiliary lemmas}
We start with the estimate on the approximation of trigonometric polynomial by simple functions.
\begin{L1}\label{pierwszy} Let $s_d,N_d\in\mathbb{N}$ and $f$ is a trigonometric polynomial. Assume that the degree with respect to the last variable of the polynomial $f$ is less than or equal to $s_d$ ($\degz{f}\leq s_d$). For a function $\tilde{f}$ given by the formula  
\[\tilde{f}(y', z) = \sum_{j=0}^{N_d-1} \chi_{_{I_j}}(z) f\left(y' , \frac{j}{N_d}\right)\qquad \forall y'\in \mathbb{T}^{d-1};\; z\in\mathbb{T},\] 
where $I_j=[\frac {j}{N_d}; \frac{j+1}{N_d}]$, we have 
\[ \left|\|f\|_{L^1(\mathbb{T}^d)}-\|\tilde{f}\|_{L^1(\mathbb{T}^d)}\right|\leq \frac {s_d}{N_d} \|f\|_{L^1(\mathbb{T}^d)}.\]
\end{L1}
\begin{proof}
We can estimate the difference of $L^1$ norms of the functions using the norm of the partial derivative of $f$.
\begin{equation*}
 \begin{split}
  \left|\|f\|_{L^1(\mathbb{T}^d)}-\|\tilde{f}\|_{L^1(\mathbb{T}^d)}\right|&=\left|\int_{\mathbb{T}^d}|f(y',z)|dy'dz-\int_{\mathbb{T}^d}|\tilde{f}(y',z)|dy'dz\right|
  \\&=\left|\sum_{j=0}^{N_d-1}\int_{\mathbb{T}^{d-1}}\int_{I_j} |f(y',z)|-|f(y',\frac{ j}{N_d})|dzdy'\right|\\
  &\leq\sum_{j=0}^{N_d-1} \int_{\mathbb{T}^{d-1}}\int_{I_j} |f(x) - f(\frac{ j}{N_d})|dzdy'
  \\&\leq \sum_{j=0}^{N_d-1}\int_{\mathbb{T}^{d-1}} \int_{I_j} \int_{0}^{z - \frac{j}{N_d}}|\frac{\partial}{\partial z} (y',z+y)|dydzdy'\\
  &\leq \sum_{j=0}^{N_d-1} \int_{I_j}\int_{\mathbb{T}^{d-1}} \int_{0}^{\frac{1}{N_d}}|\frac{\partial}{\partial z}(y',z+y)| dydy'dz\leq \frac{\|f'\|_{L^1(\mathbb{T}^d)}}{N_d}.
  \end{split}
\end{equation*}
Due to Bernstein inequality (see eg. \cite{zygmunt}) we get
\[\left|\|f\|_{L^1(\mathbb{T}^d)}-\|\tilde{f}\|_{L^1(\mathbb{T}^d)}\right|\leq  \frac{\|\frac{\partial}{\partial z} f\|_{L^1(\mathbb{T}^d)}}{N_d}\leq \frac{\degz{f}}{N}\|f\|_{L^1(\mathbb{T}^d)}\leq\frac{s_d}{N_d}\|f\|_{L^1(\mathbb{T}^d)}.\]
\end{proof}

\begin{L1}\label{drugi} For trigonometric polynomials $f_{l_d}^d, f_{l_d+1}^d \ldots, f_{k_d}^d\in L^{1}(\mathbb{T}^d)$, $-N_d<l_d<k_d\leq N_d$ and 
 \[w_{d+1}(y',y,z)= \sum_{j=l_d}^{k_d}e^{2 \pi i j y_d} f_j(y',z),\]
 following estimates are satisfied
 \begin{equation*}
  \|w_{d+1}\|_{L^1(\mathbb{T}^{d+1})}\leq \sum_{j=l_d}^{k_d}\|f_j\|_{L^1(\mathbb{T})} \leq  |k-l|  \|w_{d+1}\|_{L^1(\mathbb{T}^{d+1})} .
 \end{equation*}
\end{L1}

\begin{proof}
Left hand side of the inequality is just a triangle inequality. To get the right hand side of inequality we just observe that following inequalities are satisfied and add them up.
\[ \int_{\mathbb{T}^{d}} |f_j|dy'dz=\int_{\mathbb{T}^{d}} \left|\int_{\mathbb{T}} e^{-2 \pi i j y_d} w_{d+1}dy_d\right| dy'dz \int_{\mathbb{T}^{d}}\leq \|w_{d+1}\|_{L^1(\mathbb{T}^{d+1})}.\]
\end{proof}

\begin{L1}\label{trzeci} Assume that trigonometric polynomials $f_{l_d}^d, f_{l_{d}+1}^d \ldots, f_{k_d}^d\in L^{1}(\mathbb{T}^d)$ satisfy $\degz{f_j^d}\leq s_d$. We define functions  
 \[w_{d}(y',z):= \sum_{j=l_d}^{k_d}e^{2 \pi i N_d j z} f_j(y',z),\]
 \[w_{d+1}(y',y_d, z):= \sum_{j=l_d}^{k_d}e^{2 \pi i j y_d} f_j(y',z).\]
 This pair of functions satisfies following estimates
 \begin{equation*}
\left(1-2 |k_d-l_d| \frac{s_d}{N_d}\right)\|w_{d+1}\|_{L^1(\mathbb{T}^{d+1})}\leq\|w_d\|_{L^1(\mathbb{T}^d)}\leq \left(1+2 |k_d-l_d| \frac{s_d}{N_d}\right)\|w_{d+1}\|_{L^1(\mathbb{T}^{d+1})}.
\end{equation*}
\end{L1}

\begin{proof}
Let us define functions
\begin{equation*}
\begin{aligned}
\tilde{w}_{d+1}(y',y_d,z)&:= \sum_{j=l_d}^{k_d} e^{i j y_d} \tilde{f}_j^d(y',z),\\
\bar{w}_{d}(y',z)&:= \sum_{j=l_d}^{k_d} e^{i j N_d z} \tilde{f}_j^d (y',z).
\end{aligned} 
\end{equation*}
Using triangle inequality we get
\begin{equation}\label{oszacowanietrzeci}
\begin{split}
\left|\|w_d\|_{L^1(\mathbb{T}^d)}-\|w_{d+1}\|_{L^1(\mathbb{T}^{d+1})}\right|&\leq\left|\|w_d\|_{L^1(\mathbb{T}^d)}-\|\bar{w}_d\|_{L^1(\mathbb{T}^{d})}\right|
\\&\;+\left|\|\tilde{w}_{d+1}\|_{L^1(\mathbb{T}^{d+1})}-\|\bar{w}_d\|_{L^1(\mathbb{T}^{d})}\right|
\\&\;+\left|\|w_{d+1}\|_{L^1(\mathbb{T}^{d+1})}-\|\tilde{w}_{d+1}\|_{L^1(\mathbb{T}^{d+1})}\right|.
\end{split}
\end{equation}
Once again using the triangle inequality we obtain
\[
\left|\|w_{d}\|_{L^1(\mathbb{T}^d)}-\|\bar{w}_{d}\|_{L^1(\mathbb{T}^{d})}\right|+\left|\|w_{d+1}\|_{L^1(\mathbb{T}^{d+1})}-\|\tilde{w}_{d+1}\|_{L^1(\mathbb{T}^{d+1})}\right|\leq 2 \sum_{j=l}^k\|f_j^d - \tilde{f}_j^d\|_{L^1(\mathbb{T}^{d})}.
\]
Then the definition of the function $\tilde{f}_j^d$ leads to estimates from Lemma \ref{pierwszy} and Lemma \ref{drugi}
\begin{equation}\label{nierownosctrzeci}
\left|\|w_d\|_{L^1(\mathbb{T}^d)}-\|\bar{w}_d\|_{L^1(\mathbb{T}^{d})}\right|+\left|\|w_{d+1}\|_{L^1(\mathbb{T}^{d+1})}-\|\tilde{w}_{d+1}\|_{L^1(\mathbb{T}^{d+1})}\right|\leq   2 |k_d-l_d| \frac{s_d}{N_d}\|w_{d+1}\|_{L^1(\mathbb{T}^{d+1})}.
\end{equation}
Now we pass to estimate of the second term of the right hand side of the inequality \eqref{oszacowanietrzeci}. We know that the function  $\tilde{f}_j^d(y',\cdot)$ is a constant on a interval $I_k$ for every $j\in\{1\ldots,N_d-1\}$ and every $y'\in\mathbb{T}^{d-1}$. We denote this value by $h_j(k, y')$. This property is crucial in the following calculations.
\begin{equation*}
\begin{split}
\|\bar{w}_d\|_{L^1(\mathbb{T}^d)}&= \int_{\mathbb{T}^d} |\bar{w}(y',z)|dy'dz
\\&=\int_{\mathbb{T}^{d-1}} \sum_{k=0}^{N-1} \int_{I_k}\left|\sum_{j=l_d}^{k_d} e^{2 \pi i j N z} \bar{f}_j (y',z)\right| dzdy'
\\&=\int_{\mathbb{T}^{d-1}} \sum_{k=0}^{N-1} \int_{I_k}\left|\sum_{j=l_d}^{k_d} e^{2 \pi i j N z} h_j (k,y')\right| dzdy'
\\&=\frac{1}{N} \int_{\mathbb{T}^{d-1}} \sum_{k=0}^{N-1} \int_{\mathbb{T}}\left|\sum_{j=l_d}^{k_d} e^{2 \pi i j y_{d}} h_j(k,y')\right| dy_{d}dy'
\\&= \int_{\mathbb{T}^{d-1}}  \int_{\mathbb{T}}\sum_{k=0}^{N-1} \frac{1}{N}\left|\sum_{j=l_d}^{k_d} e^{2 \pi i j y_{d}} h_j(k,y')\right|dy_{d}dy'
\\&= \int_{\mathbb{T}^{d+1}} \left|\sum_{j=l_d}^{k_d} e^{2 \pi i j y_{d}} \bar{f}_j^d (y',z)\right|dz dy_{d}dy'
\\&= \int_{\mathbb{T}^{d+1}}|\tilde{w}_{d+1}|= \|\tilde{w}_{d+1}\|_{L^1(\mathbb{T}^{d+1})}.
\end{split}
 \end{equation*}
We have obtain
\[\|\bar{w}_d\|_{L^1(\mathbb{T}^d)}=\|\tilde{w}_{d+1}\|_{L^1(\mathbb{T}^{d+1})}\]
Above equality together with \eqref{oszacowanietrzeci} and \eqref{nierownosctrzeci} gives us 
\begin{equation*}
 \left|\|w_d\|_{L^1(\mathbb{T}^d)}-\|w_{d+1}\|_{L^1(\mathbb{T}^{d+1})}\right|\leq 2 |k_d-l_d| \frac{s_d}{N_d}\|w_{d+1}\|_{L^1(\mathbb{T}^{d+1})}
\end{equation*}
which is equivalent to the inequality from the statement of the lemma.
\end{proof}
Now using above lemmas we can prove the main theorem.
\subsection*{Proof of the main Theorem}
\begin{proof}
 Let us take polynomial $f\in L^1_F(\mathbb{T}^{\mathbb{N}})$ which depends only on first $n$ variables. Then from the definition of operator $T$ we have 
 \begin{equation*}
  Tf(z)= \sum_{(\lambda_1,\ldots,\lambda_n)\in F} \widehat{f}(\lambda_1,\ldots, \lambda_n) e^{2 \pi i \left(\sum_{i=1}^{n} \tau_j \lambda_j\right) z}\qquad \forall z\mathbb{T}, 
 \end{equation*}
which we can rewrite in the form
\begin{equation*}
  w_1(z):= Tf(z)=\sum_{j=-a_n}^{a_n} e^{2\pi i j \tau_n z} g^1_j(z), 
\end{equation*}
where $g_j$ are suitable polynomials such that $\degz{g^1_j}\leq\sum_{j=1}^{n-1} a_j |\tau_j|$. By the Lemma \ref{trzeci} we get 
\begin{equation*}
\left(1-4 a_n \frac{\sum_{j=1}^{n-1} a_j|\tau_j|}{|\tau_n|}\right)\|w_{2}\|_{L^1(\mathbb{T}^{2})}\leq\|w_1\|_{L^1(\mathbb{T})}\leq \left(1+4 a_n \frac{\sum_{i=1}^{n-1} a_j|\tau_j|}{|\tau_n|}\right)\|w_{2}\|_{L^1(\mathbb{T}^{2})},
\end{equation*}
where
\begin{equation*}
 w_2(y_1,z)=\sum_{(\lambda_1,\ldots,\lambda_n)\in F} \widehat{f}(\lambda_1,\ldots, \lambda_n) e^{2\pi i \lambda_n y_1}e^{2 \pi i \left(\sum_{i=1}^{n}\tau_j \lambda_j\right) z}. 
\end{equation*}
Analogously as in the case $d=1$ we proceed with $d>1$. We obtain trigonometric polynomials of the form
\begin{equation*}
 w_d(y',z):=\sum_{(\lambda_1,\ldots,\lambda_n)\in F} \widehat{f}(\lambda_1,\ldots, \lambda_n) e^{2\pi i \langle (\lambda_n,\ldots, \lambda_{n-d+1}), y'\rangle }e^{2 \pi i \left(\sum_{i=1}^{n-d}\tau_j \lambda_j\right) z}\qquad  \forall z\in\mathbb{T} \;\forall y' \in\mathbb{T}^{d-1}.
\end{equation*}
We can rewrite them in following way
\begin{equation*}
 w_{d}(y',z)=\sum_{j=-a_d}^{a_d} e^{2\pi i j \tau_d z} g^d_j(y',z)
\end{equation*}
with polynomials $g^d_j$ satisfying $\degz{g^d_j}\leq\sum_{j=1}^{n-d} a_j |\tau_j|$. By the Lemma \ref{trzeci} we have
\begin{equation*}
\left(1-K(d)\right)\|w_{d+1}\|_{L^1(\mathbb{T}^{d+1})}\leq\|w_d\|_{L^1(\mathbb{T}^d)} \leq \left(1+K(d)\right) 
\|w_{d+1}\|_{L^1(\mathbb{T}^{d+1})},
\end{equation*}
where constant $K(d)$ is given by the formula
\[
K(d)= 4 a_{n-d+1} \frac{\sum_{j=1}^{n-d} a_j|\tau_j|}{|\tau_{n-d+1}|}.
\]
Combining above inequalities for $d=1,...,n-1$ we get
%
%
%
%
\[
 \Pi_{j=1}^{n-1} (1-K(j))\|w_{n}\|_{L^1(\mathbb{T}^{n})}\leq\|Tf\|_{L^1(\mathbb{T})}\leq \Pi_{j=1}^{n-1} (1+K(j))\|w_{n}\|_{L^1(\mathbb{T}^{n})},
\]
Let us observe that the function $w_n$ is equal to the function $f$ up to a permutation of the variables. Hence $\|w_{n}\|_{L^1(\mathbb{T}^{n})}=\|f\|_{L^1(\mathbb{T}^{n})}$. Since $a_j, \tau_j$ satisfy $\eqref{wlasnosci}$ we have
\[ 
\sum_{d=1}^{\infty} K(d)=\sum_{d=1}^{\infty}4 a_{n-d+1} \frac{\sum_{j=1}^{n-d} a_j|\tau_j|}{|\tau_{n-d+1}|}\leq C \sum_{d=1}^{\infty} a_d (a_{d-1}+1)\frac{ \tau_{d-1}}{\tau_d}\leq\infty.
\]
\end{proof}
Hence there exist such a constant $K$ that following inequality is satisfied
\begin{equation}
  K^{-1}\|f\|_{L^1(\mathbb{T}^{n})}\leq\|Tf\|_{L^1(\mathbb{T})}\leq K\|f\|_{L^1(\mathbb{T}^{n})},
\end{equation}

\nocite{*}
\bibliographystyle{plain}

\end{document}